\documentclass[english,11pt]{amsart}
\usepackage{amsopn,amsthm,amsfonts,amsmath,amssymb,amscd}
\usepackage{babel}
\usepackage{amssymb}
\usepackage{graphicx}

\theoremstyle{definition}
\newtheorem{Theorem}{Theorem}[section]
\newtheorem{Lemma}[Theorem]{Lemma}
\newtheorem{Proposition}[Theorem]{Proposition}
\newtheorem{Corollary}[Theorem]{Corollary}
\newtheorem{Example}[Theorem]{Example}
\newtheorem{Remark}[Theorem]{Remark}

\newtheorem*{Question}{Question}

\newenvironment{Proof*}{{\it Proof.}}

\newcommand{\FF}{\mathbb{F}}

\newcommand{\QQ}{\mathbb{Q}}

\newcommand{\DEF}[1]{\emph{#1}}

\newcommand{\codim}{{\rm codim}}
\newcommand{\rk}{\mathop{\rm rk}}
\newcommand{\sh}[1]{\mathop{\rm sh}(#1)}

\newcommand{\DD}[1]{{\mathcal D}(#1)}
\newcommand{\DDi}[1]{{\mathcal D}^{-1}(#1)}

\newcommand{\n}{{\mathcal N}}
\newcommand{\cN}{{\mathcal N}}

\newcommand{\nb}{{\mathcal N}_B}
\newcommand{\nc}{{\mathcal N}_C}

\newcommand{\partition}[1]{{\mathcal P}(#1)}
\newcommand{\OO}{{\mathcal O}}

\newcommand{\rpt}[2]{r(#1,#2)}

\newcommand{\gl}[2]{\mathop{GL}_{#1}(#2)}

\begin{document}

\title{On the nilpotent commutator of a nilpotent matrix}
\author{Polona Oblak}
\date{\today}

\address{P.~Oblak: Faculty of Computer and Information Science,
Tr\v za\v ska 25, SI-1000 Ljubljana, Slovenia; e-mail: polona.oblak@fri.uni-lj.si}

\begin{abstract} 
We study the structure of the nilpotent commutator $\nb$ of a nilpotent matrix $B$. We show that $\nb$ intersects all
nilpotent orbits for conjugation if and only if $B$ is a square--zero matrix. We describe nonempty intersections  of $\nb$
with nilpotent orbits in the case the $n \times n$ matrix $B$  has rank $n-2$. Moreover, we give some results on the maximal
nilpotent orbit that $\nb$  intersects nontrivially.
\end{abstract}


\maketitle

\section{Introduction}

\bigskip

We denote by $M_n(\FF)$ the algebra of all $n\times n$ matrices over an algebraically closed field $\FF$ of characteristic
 0 and by 
$\cN_n(\FF)$ the variety of all nilpotent matrices in $M_n(\FF)$. 
Let $B\in \n_n(\FF)$ and suppose that its Jordan canonical form is given by a partition
$\underline{\lambda} \in \partition{n}$. We denote by $\nb$ the nilpotent commutator of $B$, which is the set
of all nilpotent matrices $A$ such that $AB=BA$. 
Moreover, let us denote by $\OO_B=\OO_{\underline{\lambda}}$ the orbit of $B$ under the conjugated action
of $\gl{n}{\FF}$ on $\n_n(\FF)$, i.e. the set of all nilpotent matrices with their Jordan 
canonical form given by partition $\underline{\lambda}$.

Recently, the structure of the variety of commuting nilpotent matrices has been widely studied 
(see e.g. \cite{baranovsky01,basili03,BasIar08,BaIaKh10,KoOb09}). In this paper we  investigate further which 
intersections $\nb\cap \OO_{\underline{\mu}}$ are nonempty. The answer to this question could be considered as 
a generalization of the Gerstenhaber--Hesselink Theorem on the partial order of nilpotent orbits \cite{cmc}.
 In the first part of the paper, we give anwers  for matrices 
$B$ with extremal kernel, and in the second part we give some results on the maximal partition $\underline{\mu}$, such that
the  intersection $\nb\cap \OO_{\underline{\mu}}$ is nonempty for a given $B$. 

\smallskip

In the first part of this paper (Sections 2 and 3) we are interested in  describing pairs of partitions that are  the
Jordan canonical forms of two commuting nilpotent  matrices. In Section 2, Theorem \ref{thm:2a1b}, we prove that 
the nilpotent commutator $\nb$ 
intersects every nilpotent orbit  $\OO_{\underline{\lambda}}$
if and only if $B$ is a square zero-matrix.
In Section 3 we investigate the nilpotent commutator of a nilpotent matrix having the dimension of its kernel
equal to 2. In Theorem \ref{thm:2-2}, we prove that the only pairs of distinct Jordan canonical forms of two commuting nilpotent  $n \times n$ matrices, both having exactly 2 parts,
are of the form $\left(\left(\frac{n}{2},\frac{n}{2}\right),\left(\frac{n}{2}+1,\frac{n}{2}-1\right)\right)$
where $n$ is even.
Next, we give some additional sufficient  and some necessary
conditions for partitions to be Jordan canonical forms of matrices 
in the nilpotent commutator of a nilpotent matrix, having the dimension of its kernel equal to $2$. (See
Theorems \ref{thm:antidiagonal}, \ref{thm:nilorderm+1} and Propositions \ref{thm:ind2} and also \ref{thm:ind1}.)

Some of the results in Sections 2 and 3  were already proved in \cite{PhD}. Note that recently, Britnell and
Wildon \cite{BriWil10}  proved  similar results for matrices over finite fields.

\medskip

Let us recall some definitions and notations we  use in the paper.

A nonincreasing sequence of positive integers $\underline{\mu}=(\mu_1,\mu_2,\ldots,\mu_s)$, such that
their sum is equal to $n$, is called
a \DEF{partition of an integer $n$}. It is sometimes convenient to write the partition $\underline{\mu}$ also as $(\mu_1,\mu_2,\ldots,\mu_s)=(m_1^{r_1},m_2^{r_2},\ldots,m_l^{r_l})$, 
where $\sum_{i=1}^l r_i=s$, $m_i>m_{i+1}$ and $r_i \ne 0$ for all $i$. 
By $\partition{n}$, we denote the set of all partitions of $n$.
The \DEF{conjugated partition} of a partition $\underline{\mu}$ is the partition 
$\underline{\mu^T}=(\mu_1^T, \mu_2^T,\ldots,\mu_{\mu_1}^T)$,
where $\mu_i^T=\left|\{j; \; \mu_j \geq i \}\right|$. 
It is easy to see that for each $t=1,2,\ldots,n$ there exists a uniquely defined partition 
 $\rpt{n}{t}:=(\lambda_1,\lambda_2,\ldots,\lambda_t) \in \partition{n}$, such that $\lambda_1-\lambda_t \leq 1$.
 It can be verified that $\rpt{n}{t}=\left(\left\lceil \frac{n}{t} \right\rceil^r,\left\lfloor \frac{n}{t} \right\rfloor^{t-r}\right)$. 
 By the shape of its Ferrer diagram, we call the partition $\rpt{n}{t}$ an \DEF{almost rectangular partition} of $n$.
 Moreover, we define the partial order on $\partition{n}$ with 
$(\lambda_1,\lambda_2,\ldots,\lambda_t) \leq (\mu_1,\mu_2,\ldots,\mu_s)$ if and only if 
$\sum_{i=1}^k {\lambda_i} \leq \sum_{i=1}^k {\mu_i}$ for all $k$.

Note that all eigenvalues of a nilpotent matrix are equal to 0 and thus the
Jordan canonical form of a nilpotent matrix can be described by a partition, i.e. by the decreasing sequence of 
sizes of its Jordan blocks.
If a nilpotent matrix $A$ has its Jordan canonical form given by partition $\underline{\mu}$, we write $\sh{A}=\underline{\mu}$
and call it the \DEF{shape} of matrix $A$.
For every $m$, we denote the $m \times m$ nilpotent Jordan block by $J_m$. 
By computing the lengths of the Jordan chains of $J_m^k$, $k=1,2,\ldots,m$, we observe that the Jordan canonical form of 
$J_m^k$ is given by partition $\rpt{m}{k}$. 
By $J_{\underline{\mu}}=J_{(\mu_1,\mu_2,\ldots,\mu_s)}=J_{\mu_1}\oplus J_{\mu_2}\oplus\ldots\oplus J_{\mu_s}$ we denote the 
uppertriangular matrix in its Jordan canonical form, with blocks of sizes $\mu_1 \geq \mu_2 \geq  \ldots \geq  \mu_s>0$.

\medskip

Since  $\nb$ is an irreducible variety (see Basili \cite{basili03}), there exists a unique partition 
 $\underline{\mu}$ of $n$  such that $\OO_{\underline{\mu}}\cap \nb$ is dense in $\nb$. Here, 
 $\underline{\mu}$ is the largest partition, such that the intersection 
 $\OO_{\underline{\mu}}\cap \nb$ is nonempty. Following Basili and Iarrobino \cite{BasIar08}, and Panyushev \cite{Pany2007}
 we define the map ${\mathcal D}$ on  $\partition{n}$ by 
 $\DD{\underline{\lambda}}=\underline{\mu}$. 
 
It is an interesting question (see Panyushev \cite[Problem 1]{Pany2007}) to describe $\DD{\underline{\lambda}}$ 
in terms of the partition $\underline{\lambda}$. 
Recently, some partial results to this problem were obtained. 
Basili \cite[Prop. 2.4]{basili03} showed that the number
of parts of $\DD{\underline{\lambda}}$ is equal to the smallest number $r$ such that 
$\underline{\lambda}$ is a union of $r$ almost rectangular partitions.
It was proved in \cite[Thm. 16]{oblak08} that the
first part of $\DD{\underline{\lambda}}$, where $\underline{\lambda}=(\lambda_1,\lambda_2,\ldots,\lambda_t)$, is equal to
$$ \max\limits_{1 \leq i \leq t} \left\{2(i-1)+\lambda_{i}+\lambda_{i+1}+\ldots+\lambda_{i+r}; \; 
\lambda_{i} - \lambda_{i+r}\leq 1,\ \lambda_{i-1}\ge 2\ {\rm if}\ i>1\right\}.$$ (Note that recently Ba\-si\-li and Iarrobino in \cite{BasIar_Invol}
observed the same result for an algebraically closed field $\FF$, while \cite[Thm. 16]{oblak08} holds for $\FF$ with the
characteristic 0.)

   We say that a partition $\underline{\lambda}$ is \emph{stable} if
 $\DD{\underline{\lambda}}=\underline{\lambda}$.
 Basili and Iarrobino \cite[Thm 1.12]{BasIar08} showed that $\underline{\lambda}$ is stable if and only if
its parts  differ pairwise by at least 2.
We proved in  \cite[Thm. 6]{KoOb09} that ${\mathcal D}^2={\mathcal D}$.
From these  results, we easily obtain $\DD{\underline{\lambda}}$ if it has at most two parts (see \cite[Thm. 7]{KoOb09}). 
 Until now, not much is known about $\DD{\underline{\lambda}}$ if it has more than two parts.  
In this paper, Theorem \ref{thm:dif2} characterizes  partitions $\underline{\lambda}$, such that $\DD{\underline{\lambda}}$ 
has parts that differ exactly by two. 
 In the remainder of Section 4, we examine partitions in $\DDi{\underline{\mu}}$ for certain families of partitions $\underline{\mu}$.

\bigskip
\bigskip

\section{Nilpotent commutator of a square-zero matrix}

\bigskip

We say that $B$ is a \DEF{square--zero matrix} if $B^2=0$. The Jordan canonical form of a square--zero matrix 
is given by a partition, such that all its parts are at most 2. 
By \cite[Prop.~2.4]{basili03} or \cite[Thm.~16]{oblak08}, we have that for such partition,
$\DD{(2^a,1^b)}=(2a+b)$. In the main result of this section, Theorem \ref{thm:2a1b}, we show even more: 
for an $n \times n$ square--zero matrix $B$, its nilpotent commutator $\nb$ intersects every nilpotent orbit,
i.e.~for every partition $\underline{\mu} \in \partition{n}$ there exists a nilpotent
  matrix $A$, commuting with $B$, such that $\sh{A}=\underline{\mu}$. 

\smallskip

 By $\partition{\nb}$ we denote the set of all partitions that are Jordan canonical forms of matrices in $\nb$.
 Thus, in Theorem \ref{thm:2a1b} we show that $\partition{\nb}=\partition{n}$ for every 
 $n \times n$ square--zero matrix $B$, and moreover, we show that $\partition{\nb} \subsetneqq \partition{n}$ for all 
 $n \times n$ matrices, such that $B^2\ne 0$ and $n \geq 4$.
 
 \bigskip

First, we state next Proposition, that is easy to prove, and then prove two technical lemmas
that will simplify the proof of Theorem \ref{thm:2a1b}.

\bigskip

\begin{Proposition}{\cite[Prop. 1]{oblak08}}\label{prop:AR}
A pair of partitions $((n),\underline{\mu})$
is a pair of Jordan canonical forms of two commuting $n \times n$ nilpotent matrices if and only if 
$\underline{\mu}$ is an almost rectangular
partition of $n$.
\hfill$\blacksquare$
\end{Proposition}

\bigskip

\begin{Lemma}\label{thm:3}
 If $B$ is an $n \times n$ matrix, $n \geq 4$, such that $B^2 \ne 0$, then $\partition{\nb} \subsetneqq \partition{n}$.
\end{Lemma}

\medskip
\begin{proof}
 We will show that for an arbitrary partition $\underline{\lambda} \in \partition{n}$, $\sh{B}=\underline{\lambda}$, there exists
  $\underline{\mu} \in \partition{n}$,
 such that for every $n \times n$ matrix $A$ with $\sh{A}=\underline{\mu}$, matrices $A$ and $B$ do not commute.
 
By Proposition \ref{prop:AR}, we have that if $\underline{\lambda}$
 is not an almost rectangular partition, then $(n) \notin \partition{\nb}$. 
 Suppose now $\underline{\lambda}$ is
 almost rectangular and $B^2 \ne 0$. We assume that there exists $A$, such that $\sh{A}=(n-1,1)$ and that $A$ 
 commutes with 
 $B$. We may take that $A$ is in its Jordan canonical form (otherwise, substitute $A$ with $PAP^{-1}$ and
 $B$ with $PBP^{-1}$ for  a suitable invertible matrix $P$). Then, $B$ is of the form
 $$\left[
  \begin{matrix}
   T & b\\
   c^T & 0
  \end{matrix}
 \right],$$
 where $T$ is an $(n-1) \times (n-1)$ upper triangular Toeplitz matrix and $b,c$ column vectors, such that
 $b^T=(b_1,0,\ldots,0)$ and $c^T=(0,\ldots,0,c_1)$. Let us define the matrix $B'=T \oplus 0$
  and note that $B^k=B'^k$ for all $k\geq 3$. Since $B$ is not a square-zero matrix, but $\underline{\lambda}$
 is almost rectangular, it follows that $\rk{B}> 2$. In this case, $\rk{B}=\rk{B'}$ and therefore
 $\underline{\lambda}=\sh{B}=\sh{B'}=(\rpt{n-1}{t},1)$. Since $\underline{\lambda}$
is almost rectangular, it follows that $\rpt{n-1}{t}=(2^a,1^b)$ and therefore is $B$ a square-zero matrix. This contradicts the 
assumption   and  finishes the proof that $\partition{\nb} \subsetneqq \partition{n}$.
\end{proof}

\bigskip

\begin{Lemma}\label{thm:odd}
 If $\sh{B}=(\lambda_1,\lambda_2)\in \partition{n}$, then 
 $(2^a,1^{n-2a}) \in \partition{\nb}$ for all $0 \leq a \leq \lfloor \frac{n}{2}\rfloor$.
\end{Lemma}

\smallskip

\begin{Proof*}
Suppose first that $n$ is odd. We treat the cases where $\lambda_1$ is even or $\lambda_1$ is odd, separately.
Firstly, let $\lambda_1$ be even. If $a \leq \frac{\lambda_1}{2}$, then 
 $\sh{J_{\lambda_1}^{\lambda_1-a}}=\rpt{\lambda_1}{\lambda_1-a}=(2^a,1^{\lambda_1-2a})$ and thus $\sh{J_{\lambda_1}^{\lambda_1-a} \oplus J_{\lambda_2}^{\lambda_2}}=(2^a,1^{n-2a})$.
 Obviously, $J_{\lambda_1}^{k_1}\oplus J_{\lambda_2}^{k_2}$ commutes with $J_{\lambda_1} \oplus J_{\lambda_2}$ for any positive integers $k_1$ and $k_2$.
 Otherwise, if $ \frac{\lambda_1}{2}< a \leq \frac{n-1}{2}$, then 
  $\sh{J_{\lambda_2}^{n-\frac{\lambda_1}{2}-a}}=\rpt{\lambda_2}{n-\frac{\lambda_1}{2}-a}=(2^{a-\frac{\lambda_1}{2}},1^{n-2a})$ and thus
 $\sh{J_{\lambda_1}^{\frac{\lambda_1}{2}} \oplus J_{\lambda_2}^{n-\frac{\lambda_1}{2}-a}}=(2^{a},1^{n-2a})$.
Hence,  $(2^a,1^{n-2a}) \in \partition{\nb}$
 for all $0 \leq a \leq \frac{n-1}{2}$.
 Similarly, we prove the theorem in the case $\lambda_2$ being even.

If $n$ is even, we treat  the case  $\lambda_1$ and $\lambda_2$ being even similarly as before. 
 In the case when $\lambda_1$ and $\lambda_2$ are both odd, we must consider several cases:
 \begin{itemize}
 \item  If $0 \leq a \leq \frac{\lambda_1-1}{2}$, then it is easy to see that
 $\sh{J_{\lambda_1}^{\lambda_1-a} \oplus J_{\lambda_2}^{\lambda_2}}=(2^a,1^{n-2a})$ and if $\frac{\lambda_1-1}{2} < a < \frac{n}{2}$, then
 $\sh{J_{\lambda_1}^{\frac{\lambda_1+1}{2}} \oplus J_{\lambda_2}^{\frac{\lambda_1-1}{2}+\lambda_2-a}}=(2^a,1^{n-2a})$. 
 \item If $\lambda_1$ and $\lambda_2$ are both odd and $\lambda_1=\lambda_2=a=\frac{n}{2}$, then write
 $A=\left[\begin{matrix}
  0 & I\\
  0 & 0
  \end{matrix}\right] \in \nb$, where $I \in {\mathcal M}_{a}(\FF)$.
 Since $\rk{A}=a$ and $A^2=0$, it follows that $\sh{A}=(2^a,1^{n-2a})=(2^{\frac{n}{2}})$.
 \item Suppose now that $\lambda_1$ and $\lambda_2$ are both odd, $\lambda_1> \lambda_2$, and $a=\frac{n}{2}$. 
 We write $\lambda_i=2k_i+1$, for $i=1,2$,  
 and define matrices 
 $A_{12}=\left[\begin{matrix}
  J_{\lambda_2}^{k_2}\\
  0
  \end{matrix}\right] \in {\mathcal M}_{\lambda_1 \times \lambda_2}(\FF)$, where $0 \in {\mathcal M}_{(\lambda_1-\lambda_2) \times \lambda_2}(\FF)$ 
 and
  $A_{21}=\left[\begin{matrix}
  0 & -J_{\lambda_2}^{k_2}
  \end{matrix}\right] \in {\mathcal M}_{\lambda_2 \times \lambda_1}(\FF)$, where $0 \in {\mathcal M}_{\lambda_2 \times (\lambda_1-\lambda_2)}(\FF)$.
 Here, we define $J_{\lambda_2}^0=I$.
 Let $A=\left[\begin{matrix}
  J_{\lambda_1}^{k_1} & A_{12}\\
  A_{21} & J_{\lambda_2}^{k_2+1}
  \end{matrix}\right]$. It can be easily seen that $A \in \nb$, $A^2=0$ and 
  $\rk{A}=\frac{n}{2}$. Thus $\sh{A}=(2^{\frac{n}{2}})$. \hfill$\blacksquare$
  \end{itemize}
\end{Proof*}

\bigskip

\begin{Theorem}\label{thm:2a1b}
 Let $B$ be an $n\times n$ matrix.
 \begin{itemize}
 \item If $n \leq 3$,  then $\partition{\nb}=\partition{n}$.
 \item If $n\geq 4$,
 then $\partition{\nb} = \partition{n}$ if and only if $B$ is a square--zero matrix.
 \end{itemize}
\end{Theorem}

\medskip

\begin{proof}
The case $n\leq 3$ is clear, and for $n \geq 4$ the necessity follows by Lemma \ref{thm:3}.

To prove the sufficiency of the second claim, take an arbitrary 
$\underline{\lambda}=(\lambda_1,\lambda_2,\ldots,\lambda_t) \in \partition{n}$ and let 
$B=J_{\underline{\lambda}}$. 
The matrix $B$ can be written as a direct sum 
$B_1 \oplus B_2 \oplus \ldots \oplus B_r$, where either
\begin{enumerate}
\item[(a)] $r$ is odd and $B_j=J_{\lambda_{j}}$, where all $\lambda_j$ are odd (i.e. $\underline{\lambda}$ has an odd number of odd parts), or
\item[(b)] each $B_i$ has one of the following forms: 
 \begin{enumerate}
  \item[(i)]$B_i=J_{\lambda_{j}}$, for an even $\lambda_j$,
  \item[(ii)] $B_i=J_{\lambda_{i_1}} \oplus J_{\lambda_{i_2}}$, where $\lambda_{i_1}+\lambda_{i_2}$ is even,
  \item[(iii)] $B_i=J_{\lambda_{i_1}} \oplus J_{\lambda_{i_2}}$, where $\lambda_{i_1}+\lambda_{i_2}$ is odd,
 \end{enumerate}
 and at most one $B_i$ is of the form (iii). (Namely, if $\underline{\lambda}$ has an even number of odd parts, then all $B_i$ are of the
 forms (i) and (ii), otherwise there exists exactly one $B_i$ of the form (iii).)
\end{enumerate}

It is clear that for an odd $\lambda_j$ and an arbitrary $a$, $0 \leq a \leq \frac{\lambda_j-1}{2}$, the set
$\partition{{\mathcal N}_{J_{\lambda_j}}}$ includes all partitions of the form $\rpt{\lambda_j}{\lambda_j-a}=(2^a,1^{\lambda_j-2a})$.
In the case (a), $\lambda_{2i}+\lambda_{2i+1}$ is even, thus we use Lemma \ref{thm:odd} to see that 
$\partition{{\mathcal N}_{B_{2i} \oplus B_{2i+1}}}=\partition{\lambda_{2i}+\lambda_{2i+1}}$. Therefore, 
$(2^a,1^{n-2a}) \in \partition{\nb}$ for all $a=0,1,\ldots,\lfloor \frac{n}{2}\rfloor$.

In the case (b), note that for an even $\lambda_j$ and an arbitrary $a$, $0 \leq a \leq \frac{\lambda_j}{2}$, the set
$\partition{{\mathcal N}_{J_{\lambda_j}}}$ again includes all partitions of the form $(2^a,1^{\lambda_j-2a})$.
Thus, by Lemma \ref{thm:odd}, it follows that $(2^a,1^{n-2a}) \in \partition{\nb}$ for all 
$a=0,1,\ldots,\lfloor \frac{n}{2}\rfloor$.
\end{proof}

\bigskip

\begin{Corollary}
For every nilpotent $n\times n$  matrix $A$ and integer $k \leq \frac{n}{2}$ there exists a matrix $B$, such that 
$B^2=0$ and $\rk{B}=k$.

Moreover, if $A$ and $B$ are $n \times n$ nilpotent matrices, then for each integer $k \leq \frac{n}{2}$, there 
exist a square--zero matrix $C$, such that $\rk{C}=k$, 
and $P \in \gl{n}{\FF}$, such that $C$ commutes with $A$ and $PCP^{-1}$ commutes with $B$.  
\hfill$\blacksquare$
\end{Corollary}

\bigskip

In Theorem \ref{thm:2a1b} we proved that if $B$ is not a square-zero matrix, there always exists a partition 
$\underline{\mu}$, such that the nilpotent orbit $\OO_{\underline{\mu}}$ does not intersect the nilpotent
commutator of  matrix $B$. Moreover, for a suitable $\underline{\lambda}=\sh{B}$ there exist large families
of such $\underline{\mu}$. Let us mention  the following obstruction. (See also Propositions \ref{thm:ind2}
and \ref{thm:nilorderm+1}.)

\bigskip

\begin{Proposition}\label{thm:ind1}
Let $(\underline{\lambda},\underline{\mu})$ be a pair of Jordan canonical forms of two commuting 
  nilpotent matrices, where $\underline{\lambda}=(\lambda_1,\lambda_2,\ldots,\lambda_t)$ and
  $\underline{\mu}=(\mu_1,\mu_2,\ldots,\mu_s)$.  
  If $s\geq n-\frac{\lambda_t}{2}$, then $\mu_1 \leq 2$.
\end{Proposition}
\medskip

\begin{proof}
  Let $A\in \nb$, where $\sh{B}=\underline{\lambda}$ and $\sh{A}=\underline{\mu}$ as in the statement.
   Then, $A$ can be partitioned into blocks $A_{ij}\in {\mathcal M}_{\lambda_i\times \lambda_j}(\FF)$,
   all upper triangular and constant along diagonals. 
Since $s\geq n-\frac{\lambda_t}{2}$, we have $\rk A=n-s\leq \frac{\lambda_t}{2}$ and thus for all $i,j$,
 $\rk A_{i,j}\leq \frac{\lambda_t}{2}$. It follows that $A^2=0$ and thus $\mu_1 \leq 2$.
\end{proof}

\bigskip

\begin{Example}\label{ex:transpose}
Note that Baranovsky proved in \cite[Lemma 3]{baranovsky01} that 
$(\underline{\lambda},\underline{\lambda}^T)$ is a pair of Jordan canonical forms of two commuting nilpotent
matrices.

In the case $\underline{\lambda}=(\lambda_1,\lambda_2)$, $\underline{\lambda}^T$ has all parts equal to at most two and thus by Theorem \ref{thm:2a1b}, $(\underline{\lambda},\underline{\mu})$
is a pair of Jordan canonical forms of two commuting nilpotent matrices for all 
$\underline{\mu}\leq \underline{\lambda}^T$. However, this is not true in general.

\smallskip

 Suppose $\sh{B}=\underline{\lambda}=(\lambda_1,\lambda_2,\ldots,\lambda_t)$, where $t\geq 3$ and
  $\lambda_t \geq 4$. Then, $\left(\underline{\lambda},(3,1^{n-3})\right)$ is {\bf not} a pair of Jordan 
  canonical forms of two commuting nilpotent matrices (see Proposition \ref{thm:ind1})
  and $(3,1^{n-3}) \leq \underline{\lambda}^T$.
  \hfill$\square$
\end{Example}

\bigskip
\bigskip

\section{Partitions with 2 parts}

\bigskip

Besides the square--zero matrices that have rather large dimension of its kernel, we are also interested in 
matrices, having its kernel of dimension at most two, i.e. matrices 
that have at most two Jordan blocks. Jordan canonical forms of matrices in the nilpotent commutator  of 
the matrix with one Jordan block are characterized in Proposition \ref{prop:AR}.

In this section, we give a characterization of pairs of Jordan canonical forms of two commuting nilpotent 
matrices, each having exactly two Jordan blocks. Namely, we will prove the following.

\bigskip

\begin{Theorem} \label{thm:2-2}
 A pair $\left((\lambda_1,\lambda_2),(\mu_1,\mu_2)\right)$ of \emph{distinct} partitions of $n$ is a pair of
  Jordan canonical forms of two nilpotent commuting matrices if and only if 
 $n$ is even and one of them is equal to $\left(\frac{n}{2},\frac{n}{2}\right)$ and the other one is equal to 
 $\left(\frac{n}{2}+1,\frac{n}{2}-1\right)$.
\end{Theorem}

\bigskip

%

 Define matrices $M=J_{\lambda_1}\oplus 0$, where $0 \in {\mathcal M}_{\lambda_2}(\FF)$,
and $N= 0 \oplus J_{\lambda_2}$, where $0 \in {\mathcal M}_{\lambda_1}(\FF)$. 
Write also $M_0=I_{\lambda_1 \times \lambda_1} \oplus 0_{\lambda_2 \times \lambda_2}$
and $N_0=0_{\lambda_1 \times \lambda_1} \oplus I_{\lambda_2 \times \lambda_2}$.
Let us write $M_i=M^i$ and $N_i=N^i$ for $i=0,1,\ldots$.

For $k=0,1,\ldots,\lambda_2-1$ let $K_k$ be an $n \times n$ matrix such that its only nonzero entries are in the 
positions $(i,\lambda_1+k+i)$,
where $i=1,2,\ldots,\lambda_2-k$, and are all equal to 1. Similarly, let us define matrices $L_l$ for $l=0,1, \ldots,\lambda_2-1$ 
such that its only nonzero entries (which are equal to 1) are in the positions $(\lambda_1+j,\lambda_1-\lambda_2+l+j)$, where 
$j=1,2,\ldots,\lambda_2-l$.
       
\bigskip
  
It is easy to see that the only nonzero products of these matrices are:
  \begin{equation}\label{array1}
     \begin{array}{rlcrl}
        M_i \cdot M_j&=M_{i+j} &\qquad \qquad&  M_i \cdot K_j&=K_{i+j}\\
        K_i \cdot L_j&=M_{\lambda_1-\lambda_2+i+j} &\qquad \qquad&   K_i \cdot N_j&=K_{i+j}\\
        L_i \cdot M_j&=L_{i+j} &\qquad \qquad& L_i \cdot K_j&=N_{\lambda_1-\lambda_2+i+j}\\
        N_i \cdot L_j&=L_{i+j} &\qquad \qquad& N_i \cdot N_j&=N_{i+j}\\
      \end{array}
    \end{equation}
where by the definition $M_j=0$ for all $j \geq \lambda_1$ and $K_{i}=L_{i}=N_{i}=0$ for $i \geq \lambda_2$.

\bigskip

From now on, let $B=J_{\lambda_1}\oplus J_{\lambda_2}$, where $\lambda_1 \geq \lambda_2 > 0$.

It is well known that nilpotent matrix $A$, commuting with $B$, is of the form
\begin{equation}\label{eq:longA}
   A=\sum_{i=1}^{\lambda_1-1}a_i M_i+\sum_{i=0}^{\lambda_2-1}b_i K_i+ \sum_{i=0}^{\lambda_2-1}c_i L_i+\sum_{i=1}^{\lambda_2-1}d_i N_i\, 
\end{equation}
where $b_0 c_0=0$ if $\lambda_1=\lambda_2$. Equivalently,
\begin{equation*}
   A=\sum_{i=\alpha}^{\lambda_1-1}a_i M_i+\sum_{i=\beta}^{\lambda_2-1}b_i K_i+ \sum_{i=\gamma}^{\lambda_2-1}c_i L_i+\sum_{i=\delta}^{\lambda_2-1}d_i N_i\, ,
\end{equation*}
where $a_{\alpha} b_{\beta} c_{\gamma} d_{\delta} \ne 0$ and if $\lambda_1=\lambda_2$, also $\beta+\gamma \geq 1$.
We define $\alpha=\lambda_1$ (resp. $\beta=\lambda_2$, $\gamma=\lambda_2$, $\delta=\lambda_2$) if $a_i=0$ for $i=1,2,\ldots,\lambda_1-1$
(resp. $b_i=0$ for $i=0,1,\ldots,\lambda_2-1$, $c_i=0$ for $i=0,1,\ldots,\lambda_2-1$, $d_i=0$ for $i=1,2\ldots,\lambda_2-1$).


\bigskip

In what follows, we will prove some lemmas that will give the proof of Theorem \ref{thm:2-2}. Lemma
\ref{thm:rk} is well known, but we give here full proof for the sake of completeness.

\bigskip

\begin{Lemma}\label{thm:rk}
 If $A$ is a nilpotent matrix, such that $\sh{A}=(\lambda_1,\lambda_2,\ldots,\lambda_t)$, then
 the kernel of matrix $A^j$ has dimension equal to $\sum_{i=1}^j \lambda_i^T$.
\end{Lemma}

\smallskip

\begin{proof}
 For a nilpotent matrix $A$ with $\sh{A}=\underline{\lambda}=(\lambda_1,\lambda_2,\ldots,\lambda_t)$ let 
 $V_{\lambda_1} \oplus V_{\lambda_2} \oplus \ldots \oplus V_{\lambda_t}$ be a decomposition of $\FF^n$ corresponding to 
 the Jordan canonical form $J_{\underline{\lambda}}=J_{\lambda_1}\oplus J_{\lambda_2}\oplus\ldots\oplus J_{\lambda_t}$ of $A$.
 For matrix $J_{\underline{\lambda}}$ it is clear that  
 $$\codim _{\ker J_{\underline{\lambda}}^{j}|_{V_{\lambda_i}}} (\ker J_{\underline{\lambda}}^{j-1}|_{V_{\lambda_i}})=\left\{ 
    \begin{array}{cl}
      1, & \text{ if } \lambda_i \geq j, \\
      0, & \text{ otherwise. }
    \end{array} \right.$$
 Therefore $\codim_{\ker{J_{\underline{\lambda}}^{j}}} (\ker{J_{\underline{\lambda}}^{j-1}})=|\{\lambda_i; \; \lambda_i \geq j\}|= \lambda_j^T$ and thus
 $\dim \ker A^j=\dim \ker J_{\underline{\lambda}}^j=
 \sum_{i=1}^j \codim _{\ker J_{\underline{\lambda}}^{i}} (\ker J_{\underline{\lambda}}^{i-1})=\sum_{i=1}^j \lambda_i^T$.
\end{proof}

\bigskip

\begin{Lemma} \label{thm:geq3}
%
Suppose that  either $\lambda_1-\lambda_2=1$ or $\lambda_1-\lambda_2 \geq 3$.
The pair $\left((\lambda_1,\lambda_2),(\mu_1,\mu_2)\right)$ is a pair 
of Jordan canonical forms of two commuting nilpotent matrices
if and only if $(\mu_1,\mu_2) = (\lambda_1,\lambda_2)$.
\end{Lemma}

   \smallskip
 
\begin{proof}
 Let $B=J_{\lambda_1}\oplus J_{\lambda_2}$ and suppose there exists $A \in \nb$, such that
  $\sh{A}=\underline{\mu}=(\mu_1,\mu_2)\ne (\lambda_1,\lambda_2)$. 
 
 Suppose first that $\lambda_1-\lambda_2 \geq 3$. Since $(\lambda_1,\lambda_2)$ is stable, we have that  
 $\mu_1 < \lambda_1$ and $\mu_2 \geq \lambda_2+1$. 
 By Lemma \ref{thm:rk},
 $\rk {A^{\lambda_2+1}} = n- 2 (\lambda_2+1)=\lambda_1-\lambda_2-2$. 
 On the other hand, since $A \in \nb$ is of the form \eqref{eq:longA} and $\lambda_1-\lambda_2 \geq 2$, it follows that
 $a_1d_1 \ne 0$.  Since $\lambda_1-\lambda_2 \geq 3$, it follows that 
 $A^{\lambda_2+1}=\sum_{i=\lambda_2+1}^{\lambda_1-1}a'_i M_i$, where $a'_{\lambda_2+1}\ne 0$. Thus, 
 $\rk{A^{\lambda_2+1}}=\lambda_1-\lambda_2-1$, which is a contradiction.
 
 If $\lambda_1-\lambda_2=1$, we have that $\mu_1 -\mu_2 \geq 3$. It follows from the previous paragraph that
 $((\mu_1,\mu_2),(\lambda_1,\lambda_2))$ is not a pair of Jordan canonical forms of two commuting nilpotent matrices.
\end{proof}

\bigskip

\begin{Lemma} \label{thm:=2}
 If $\lambda_1-\lambda_2 =2$ and $\left((\lambda_1,\lambda_2),(\mu_1,\mu_2)\right)$ is a pair 
of Jordan canonical forms of two commuting nilpotent matrices, then
 $\mu_1 = \lambda_1$ or $\mu_1=\lambda_1-1$.
%
 \end{Lemma}

   \smallskip

\begin{proof}
 To prove the lemma, it suffices to prove that $(\lambda_1-1,\lambda_1-1) \in \partition{\nb}$
 for $\sh{B}=(\lambda_1,\lambda_1-2)$. Equivalently, we
  have to prove that  $(\lambda+1,\lambda-1) \in \partition{\nc}$, where $\sh{C}=(\lambda,\lambda)$. 
 
 Define an upper triangular matrix 
 $A=\sum_{i=1}^{\lambda-1}a_i M_i+\sum_{i=0}^{\lambda-1}b_i K_i+ \sum_{i=1}^{\lambda-1}d_i N_i \in \nc$, where $a_1$, $b_0$ and $d_1$ are
 algebraically independent over $\QQ$.  
 From  \eqref{array1}, it easily follows that 
 $A^k=\sum_{i=k}^{\lambda-1}a'_i M_i+\sum_{i=k-1}^{\lambda-1}b'_i K_i+ \sum_{i=k}^{\lambda-1}d'_i N_i$.
 Thus, $\rk{A^k}= n-2k$ for $k=1,2,\ldots,\lambda-1$, $\rk{A^{\lambda}}=1$ and $A^{\lambda+1}=0$. 
 Therefore, $\sh{A}=(\lambda+1,\lambda-1)$ and thus $(\lambda+1,\lambda-1) \in \partition{\nc}$.
\end{proof}

\bigskip

Now, Theorem \ref{thm:2-2} follows from Lemmas \ref{thm:geq3} and \ref{thm:=2}.
As a Corollary of Theorem \ref{thm:2-2}, we get the partitions of maximal rank in $\partition{\nb}$.

\bigskip

\begin{Corollary} 
 Let  $\sh{B}=(\lambda_1,\lambda_2)\in \partition{n}$. The set of Jordan canonical forms of matrices $A \in \nb$ of maximal rank is equal to
  \begin{enumerate}
   \item[(a)]  $\{(n)\}$ if $\lambda_1-\lambda_2\leq 1$,
   \item[(b)]  $\{(\lambda_1,\lambda_2), (\lambda_1-1,\lambda_2+1)\}$ if $\lambda_1-\lambda_2 =2$,
   \item[(c)]  $\{(\lambda_1,\lambda_2)\}$ if $\lambda_1-\lambda_2 \geq 3$. \hfill$\blacksquare$
 \end{enumerate}
\end{Corollary}

\bigskip

Next, we add some Jordan canonical forms of matrices in the nilpotent commutator of matrix 
$J_{\lambda_1} \oplus J_{\lambda_2}$, which are not almost rectangular subpartitions of $(\lambda_1,\lambda_2)$.

\bigskip

\begin{Theorem} \label{thm:antidiagonal}
Let $\sh{B}=(\lambda_1,\lambda_2)\in \partition{n}$.
Choose integers $j,\ell$ such that $0 \leq j \leq \ell < \lambda_2$ and write $w=\lambda_1-\lambda_2+j+\ell$. 
 \begin{enumerate}
  \item[(a)] If there exists an integer $k$, such that $\lambda_2\leq kw < \lambda_1$,
   then $$\left((2k+1)^{\lambda_1-kw},(2k)^{w+\lambda_2-\lambda_1},(2k-1)^{kw-\lambda_2}\right) \in \partition{\nb}\, .$$
  \item[(b)] If there exists an integer $k$, such that $\lambda_1-\ell \leq kw <\lambda_2-j$,
   then $$\left((2k+2)^{\lambda_2-kw-j},(2k+1)^{w+j-\ell},(2k)^{kw+\ell-\lambda_2}\right) \in 
   \partition{\nb}\, .$$
  \item[(c)] Otherwise, $\rpt{n}{w} \in \partition{\nb}$.
 \end{enumerate}
\end{Theorem}

  \medskip
 
\begin{Proof*}
Given $j, \ell$ with the desired properties, let us define the matrix 
$A= b K_j+c L_\ell
$ and $w=\lambda_1-\lambda_2+j+\ell$. For such $A$, using induction on $m$, it is easy to verify that for all 
$m \geq 1$:
  \begin{equation}\label{array2}
     \begin{array}{rlcrl}
 \rk{(A^{2m})_{11}} & =\max\{\lambda_1-mw,0\} &\qquad \qquad& (A^{2m-1})_{11}&=0\\
 \rk{(A^{2m})_{22}}&=\max\{\lambda_2-mw,0\}&\qquad \qquad&(A^{2m-1})_{22}&=0 \\
 \rk{(A^{2m-1})_{12}}&=\max\{\lambda_1+\ell-mw,0\}&\qquad \qquad& (A^{2m})_{12}&=0\\
 \rk{(A^{2m-1})_{21}}&=\max\{\lambda_1+j-mw,0\}&\qquad \qquad& (A^{2m})_{21}&=0
 \end{array}
\end{equation}
 Note that $\rk{(A^{2m})_{11}} \geq \rk{(A^{2m})_{22}}$ and
 $\rk{(A^{2m-1})_{21}} \leq \rk{(A^{2m-1})_{12}}$ for all $m \geq 1$, since $j\leq \ell$.

 \begin{enumerate}
  \item[(a)] First, suppose there exists an integer $k$ such that  $\lambda_1>k w$ and $\lambda_2 \leq k w$. By 
 \eqref{array2}, we have that  $\rk{(A^{2k})_{11}}>0$
   and $\rk{(A^{2k})_{22}}=0$. Since 
   $\rk{(A^{2k+1})_{12}}=0$, it follows that $A^{2k+1}=0$.
   On the other hand $\rk(A^{2k-1})_{12} \geq \rk{(A^{2k-1})_{21}}>0$.

    If $m$ is even and $m \leq 2k-1$, then by \eqref{array2} it follows that
    $\dim \ker A^{m}= \lambda_1+\lambda_2-\rk{(A^{m})_{11}}-\rk(A^{m})_{22}=mw$.
    Similarly, if $m \leq 2k-1$ is odd, then
    $\dim \ker A^{m}=\lambda_1+\lambda_2-\rk{(A^{m})_{12}}-\rk(A^{m})_{21}=mw$.
     Since $kw < \lambda_1$, it follows that $\lambda_2-(k-1)w > \lambda_1-\lambda_2+w \geq 0$ and therefore
   \begin{eqnarray*}
    \sh{A}&=&\left(w^{2k-1},\lambda_2-(k-1)w,\lambda_1-kw\right)^T=\\
          &=&\left((2k+1)^{\lambda_1-kw},(2k)^{w+\lambda_2-\lambda_1},(2k-1)^{kw-\lambda_2}\right)\, .
   \end{eqnarray*}

  \item[(b)] If there exists an integer $k$ such that 
   $\lambda_1-\ell \leq kw < \lambda_2-j$, we proceed similarly as in (a) and observe that 
   $\rk{(A^{2k+1})_{12}}>0$, $\rk{(A^{2k+1})_{21}}= 0$,
   and $\rk(A^{2k-1})_{12} \geq \rk{(A^{2k-1})_{21}}>0$. Since 
   $\rk{(A^{2k+2})_{11}}=0$ 
   it follows by \eqref{array2} that $A^{2k+2}=0$.
   Again, similarly as in (a), we can compute that
   $\dim \ker (A^{m})=mw$ for all $1 \leq m \leq 2k$ and thus
   \begin{eqnarray*}
    \sh{A}&=&\left(w^{2k},\lambda_1+j-kw,\lambda_2-j-kw\right)^T=\\
     &=&\left((2k+2)^{\lambda_2-kw-j},(2k+1)^{w+j-\ell},(2k)^{kw+\ell-\lambda_2}\right).
   \end{eqnarray*}
   
  \item[(c)] Since there does not exist an integer $k$ such that 
   $\lambda_2 \leq kw < \lambda_1$, it follows that $\rk{(A^k)_{11}}=0$ if and only if $\rk{(A^k)_{22}}=0$.
   Similarly, since there does not exist an integer $k$ such that 
   $\lambda_1-\ell \leq kw < \lambda_2-j$, we can conclude that
   $\rk{(A^k)_{12}}=0$ if and only if $\rk{(A^k)_{21}}=0$ for all $k$. Write $n=sw+r$, where
   $0 \leq r < w$. Similarly as in (a), we conclude that 
   $\dim \ker (A^{i})=iw$ for all $1 \leq i \leq s$ and
   $A^{s+1}=0$. Therefore $\sh{A}=\left(w^s,r \right)^T=\rpt{n}{w}$.
   \hfill$\blacksquare$
 \end{enumerate}
\end{Proof*}

\bigskip

On the other hand, there exist plenty of partitions, that are not Jordan canonical forms of matrices 
in the nilpotent commutator $\nb$ for $\sh{B}=\underline{\lambda}$. (We already proved Proposition
\ref{thm:ind1}.)

\bigskip

The following lemma can be verified straightforwardly and will be used to prove Proposition \ref{thm:ind2}.

\bigskip

\begin{Lemma}\label{thm:rank}
\begin{enumerate}
\item
 If $C \in {\mathcal M}_{p \times r}(\FF)$ and $D \in {\mathcal M}_{r \times q}(\FF)$ are uppertriangular matrices,
 constant along diagonals, 
 then their product $CD$ is also of the same form and $\rk(CD)=\max \{\rk(C)+\rk(D)-r,0\}$.
 \item
  If $A=\sum_{i=\alpha}^{\lambda_1-1}a_i M_i+\sum_{i=\beta}^{\lambda_2-1}b_i K_i+\sum_{i=\gamma}^{\lambda_2-1}c_i L_i+ \sum_{i=\delta}^{\lambda_2-1}d_i N_i$, where $a_\alpha b_\beta c_\gamma d_\delta \ne 0$, 
 then $\rk{A}\leq \max\{n-\alpha-\delta,2\lambda_2-\beta-\gamma\}$.
 \hfill$\blacksquare$
 \end{enumerate}
\end{Lemma}

\bigskip

\begin{Proposition} \label{thm:ind2}
Let $\left((\lambda_1,\lambda_2),(\mu_1,\mu_2,\ldots,\mu_s)\right)$ be a pair of Jordan canonical forms
of two commuting nilpotent matrices. If $s > \lambda_1$, then  
$\mu_1 \leq \left\lceil \frac{\lambda_2}{s-\lambda_1} \right\rceil$.
\end{Proposition}

\smallskip
 
\begin{proof}
  Let $\sh{B}=(\lambda_1,\lambda_2)$ and 
  $A=\left[
\begin{matrix}
 A_{11} & A_{12}\\
 A_{21} & A_{22}
\end{matrix}
\right] \in \nb$, 
 $A_{ij}\in {\mathcal M}_{\lambda_i \times \lambda_j}(\FF)$. Note that $A_{ij}$, $1\leq i,j \leq 2$, are upper 
 triangular 
 and constant along diagonals. Let $\sh{A}=(\mu_1,\mu_2,\ldots,\mu_s)$ and denote
  $q=\rk{A}=n-s$.  It follows that $\rk{A_{ij}}\leq q$ for all $1 \leq i,j \leq 2$. By assumption, $s>\lambda_1$ and thus 
 $q < \lambda_2$.
 
  We first prove that $\rk{(A^k)_{ij}}\leq \max\{kq-(k-1)\lambda_2,0\}$ for all $k \geq 1$. 
 
 The case $k=1$ is clear and then we proceed by induction. Suppose that 
 $\rk{(A^k)_{ij}}\leq \max\{kq-(k-1)\lambda_2,0\}$ for all $1 \leq i,j \leq 2$. 
 If $kq-(k-1)\lambda_2 \geq 0$, then $(k+1)q-k \lambda_2 \geq q-\lambda_2$. Therefore, by Lemma \ref{thm:rank}, 
 \begin{eqnarray*}
  \rk{(A^{k+1})_{ij}} &\leq& \max\big\{\rk\big((A^k)_{i1}A_{1j}), \rk((A^k)_{i2}A_{2j}\big)\big\}\leq\\
     &\leq& \max\{kq-(k-1)\lambda_2+q-\lambda_1, q-\lambda_1, \\
      && \quad \quad \; kq-(k-1)\lambda_2+q-\lambda_2,q-\lambda_2\}=\\
     &\leq& \max\{(k+1)q-k \lambda_2,q-\lambda_2\}=\\
     &=&(k+1)q-k \lambda_2 \leq \max\{(k+1)q-k \lambda_2,0\}
 \end{eqnarray*}
 for all $1 \leq i,j \leq 2$. 
 If $kq-(k-1)\lambda_2 \leq 0$, then $(A^k)_{ij}=0$ for $1 \leq i,j \leq 2$ and thus $\rk{(A^{k+1})_{ij}} = 0$ for all 
 $1 \leq i,j \leq 2$.
 
 Again, using Lemma \ref{thm:rank} it follows that 
 \begin{eqnarray*}
  \rk{A^k} &\leq& \max\{\rk{A^{k}}_{11}+\rk{A^{k}}_{22},\rk{A^{k}}_{12}+\rk{A^{k}}_{21}\} \leq \\
           &\leq& 2\max\{kq-(k-1)\lambda_2,0\} = \\
           &=& \max\{2\lambda_2-2k(\lambda_2-q),0\}
 \end{eqnarray*}
 for all $k \geq 1$. 
 
 Thus for all $k \geq \frac{\lambda_2}{\lambda_2-q}$ it follows that $\rk{A^k} = 0$. Therefore, 
 $\mu_1 \leq \left\lceil \frac{\lambda_2}{\lambda_2-q} \right\rceil=\left\lceil \frac{\lambda_2}{s-\lambda_1} \right\rceil$.
\end{proof} 

\bigskip

 \begin{Proposition} \label{thm:nilorderm+1}
  Let $\sh{B}=(\lambda,\lambda)$ and $A \in \nb$. 
  If a partition $\underline{\mu}=(\mu_1,\mu_2, \ldots,\mu_s) \in \partition{\nb}$, then
  either $\underline{\mu}=(n)$ or $\mu_1 \leq \lambda+1$.
\end{Proposition}
    
    \smallskip
 
 \begin{proof}
  Write $A \in \nb$ as in \eqref{eq:longA}, $b_0c_0=0$, and assume that $A^{n-1}=0$.
   Without loss of generality suppose that $c_0=0$. 
  
  Since $A^{n-1}=A^{2\lambda-1}=b_0^{\lambda} c_1^{\lambda-1} K_{\lambda-1}$, it follows that 
  $b_0c_1=0$.
  If $b_0=0$, then $A=\sum_{i=1}^{\lambda-1}a_i M_i+\sum_{i=1}^{\lambda-1}b_i K_i+\sum_{i=1}^{\lambda-1}c_i L_i+ \sum_{i=1}^{\lambda-1}d_i N_i$
  and by \eqref{array1}, there are no summands in $A^{\lambda}$ having index less than $\lambda$. Therefore, 
  $A^{\lambda}=0$.    
  In the case $c_1=0$, the only summands of $A^{\lambda+1}$ having smaller index that the number of factors,
  are $(a_1M)^\lambda  b_0K_0=0$ and $b_0K_0 (d_1N)^\lambda =0$. Thus $A^{\lambda+1}=0$. 
\end{proof}
  
%
%
%
%

\bigskip
\bigskip

\section{Inverse image of $\DD{\underline{\lambda}}$ 
}

\bigskip

In the case, when $\DD{\underline{\lambda}}$ has at most 2 parts, the partition $\DD{\underline{\lambda}}$ can be easily characterized 
in the terms of $\underline{\lambda}$ (see \cite[Thm. 7]{KoOb09}). There is not much known about 
$\DD{\underline{\lambda}}$ if it has at least 3 parts.


\bigskip

\begin{Theorem}\label{thm:dif2}
 For a partition $\underline{\lambda}$,
 $$\DD{\underline{\lambda}}=(\mu,\mu-2,\mu-4,\ldots,\mu-2k)$$
 if and only if 
  $$\underline{\lambda}=(\mu,\mu-2,\mu-4,\ldots,\mu-2k+2,\rpt{\mu-2k}{t}),$$
 for $t=1,2,\ldots,\mu-2k$.
 
 Therefore, $\left| \DDi{\mu,\mu-2,\mu-4,\ldots,\mu-2k} \right|=\mu-2k$.
\end{Theorem}

\smallskip

\begin{proof}
 Suppose $\DD{\underline{\lambda}}=(\mu,\mu-2,\mu-4,\ldots,\mu-2k)$ and 
 obtain that $\underline{\lambda}$ is a partition of  $n=(k+1)\mu-k(k+1)$. 
 Since $\DD{\underline{\lambda}}$ has $k+1$ parts,  it follows by \cite[Thm.~2.4]{basili03}, that 
 $\underline{\lambda}$ is of  the form
 $$\underline{\lambda}=(\lambda_{1,1},\lambda_{1,2},\ldots,\lambda_{1,t_1},
               \lambda_{2,1},\lambda_{2,2},\ldots,\lambda_{2,t_2},\ldots,
               \lambda_{k+1,1},\lambda_{k+1,2},\ldots,\lambda_{k+1,t_{k+1}})\, ,$$
 where $(\lambda_{i,1},\lambda_{i,2},\ldots,\lambda_{i,t_i})$ are almost rectangular partitions and $t_i\geq 1$,
 $i=1,2,\ldots,k+1$.
 Since the first part of $\DD{\underline{\lambda}}$ is equal to $\mu$, it follows by \cite[Thm. 16]{oblak08} that
 \begin{eqnarray}\label{eq:path}
  \lambda_{1,1}+\lambda_{1,2}+\ldots+\lambda_{1,t_1} &\leq& \mu\, , \notag \\
  2t_1+\lambda_{2,1}+\lambda_{2,2}+\ldots+\lambda_{2,t_2}  &\leq& \mu\, , \notag \\
  2(t_1+t_2)+\lambda_{3,1}+\lambda_{3,2}+\ldots+\lambda_{3,t_3}  &\leq& \mu\, , \\
   &\vdots& \notag \\
  2(t_1+t_2+\ldots+t_k)+\lambda_{k+1,1}+\lambda_{k+1,2}+\ldots+\lambda_{k+1,t_{k+1}}  &\leq& \mu\, , \notag
 \end{eqnarray} 
 where at least one of the inequalities is actually an equality.
 
 By summing all inequalities, we have
  $2(kt_1+ (k-1)t_2+\ldots+2t_{k-1}+t_{k})+n\leq (k+1)\mu$.
 Since $t_i \geq 1$ for all $i$, it follows that 
  $k(k+1) \geq 2(kt_1+ (k-1)t_2+\ldots+2t_{k-1}+t_{k})\geq 2 (k+(k-1)+\ldots+2+1)=k(k+1)$, 
 and therefore, $t_i=1$ for $i=1,2,\ldots,k$ and all inequalities in \eqref{eq:path} are equalities. 
 By the last inequality in \eqref{eq:path} it follows  that
  $ \lambda_{k+1,1}+\lambda_{k+1,2}+\ldots+\lambda_{k+1,t_{k+1}}  = \mu-2k$.
 
 Now, $\underline{\lambda}$ has the form $\underline{\lambda}=(\lambda_{1}, \lambda_{2}, \ldots, \lambda_{k},
               \lambda_{k+1,1},\lambda_{k+1,2},\ldots,\lambda_{k+1,t_{k+1}})$,
 and since $\DD{\underline{\lambda}}$ has $k+1$ parts, it follows that  $\lambda_{i-1} -\lambda_i \geq 2$
 for $i=1,2,\ldots,k$. Suppose there exists $j$, $2 \leq j \leq k$, such that 
 $\lambda_j < \lambda_1-2(j-1)$ and let $j$ be minimal such. Thus, for $i=1,2,\ldots,j-1$ we have  
 $\lambda_i=\lambda_1-2(i-1)$ and for $i=j,j+1,\ldots,k$, we have $\lambda_i\leq \lambda_1-2i+1$. Using
 the above equalities, we  now obtain
 \begin{eqnarray*}
 k \lambda_1 -k(k-1) &\leq& k\mu-k(k-1) =n-(\mu-2k)= \\
     &=&  n-(\lambda_{k+1,1}+\lambda_{k+1,2}+\ldots+\lambda_{k+1,t_{k+1}})=\\
     &=& \sum_{i=1}^k \lambda_i = \sum_{i=1}^{j-1} \lambda_i + \sum_{i=j}^k \lambda_i\leq\\
     &\leq&k\lambda_1-k(k-1)-(k-j+1)
 \end{eqnarray*}
 and therefore  $j \geq k+1$, which contradicts the existence of $j$, $2 \leq j \leq k$, such that 
 $\lambda_j < \lambda_1-2(j-1)$. Thus, $\lambda_i=\lambda_1-2(i-1)$ for $i=1,2,\ldots,k$ and thus 
 $(\lambda_{k+1,1},\lambda_{k+1,2},\ldots,\lambda_{k+1,t_{k+1}})$ is an almost rectangular partition of 
 $\mu-2k$.
\end{proof}

\bigskip

\begin{Remark}
Note that Theorem \ref{thm:dif2} does not hold if the parts of  
$\DD{\underline{\lambda}}$ differ for at least 3. For example, $\DD{(3,1,1)}=(4,1)$.
\end{Remark}

\bigskip

As mentioned in the begining of this section, $\DD{\underline{\lambda}}$ is known when $\underline{\lambda}$ has at most
 2 parts. (See \cite[Thm. 7]{KoOb09}). Here, we describe the preimage of ${\mathcal D}$ for certain partitions and give a conjecture
 on the size of  $\DDi{\underline{\mu}}$ in the case $\underline{\mu}$ has two parts. We will need the following 
 lemma.
 
\bigskip

\begin{Lemma}\label{thm:lemma1}
 If $\DD{\underline{\lambda}}=(\mu,\mu-r)$, where $2 \leq r < \mu$, then the partition $\underline{\lambda}$ is of the form
 $\underline{\lambda}=(\lambda_1,\lambda_2,\ldots,\lambda_s,\lambda_{s+1},\ldots,\lambda_t)$,
 where 
 \begin{itemize}
  \item $\lambda_1-\lambda_s \leq 1$, $\lambda_{s+1}-\lambda_t \leq 1$,
  \item $\lambda_1-\lambda_t \geq 2$,
  \item $s\leq \frac{r}{2}$.
 \end{itemize} 
\end{Lemma}
 
 \smallskip

\begin{proof}
 If $\DD{\underline{\lambda}}=(\mu,\mu-r)$, then  $\underline{\lambda}$ is of the form  
 $\underline{\lambda}=(\lambda_1,\lambda_2,\ldots,\lambda_s,\lambda_{s+1},\ldots,\lambda_t)$, where $\lambda_1-\lambda_s \leq 1$,
 $\lambda_{s+1}-\lambda_t \leq 1$ and $\lambda_1-\lambda_t \geq 2$. (See Basili \cite[Prop. 2.4]{basili03}.) 
 Since the first part of $\DD{\underline{\lambda}}$ is equal to $\mu$, it follows by \cite[Thm. 16]{oblak08} that 
 \begin{eqnarray}
  \lambda_1+\lambda_2+\ldots+\lambda_s &\leq& \mu \label{eq:11}\\
  2s+\lambda_{s+1}+\lambda_{s+2}+\ldots+\lambda_t &\leq& \mu \, .\label{eq:12}
 \end{eqnarray} 
 Thus, $2\mu-r=\lambda_1+\lambda_2+\ldots+\lambda_t\leq 2\mu-2s$ and 
 therefore $s\leq \frac{r}{2}$.
\end{proof}

\bigskip

\begin{Proposition}
 For a partition $\underline{\lambda}=(\lambda_1,\lambda_2,\ldots,\lambda_t) \in \partition{n}$ and 
 $n \geq 4$, it follows that 
 $$\DD{\underline{\lambda}}=(n-1,1)$$
 if and only if 
  $\lambda_1-\lambda_t\geq 2$ and either $\underline{\lambda}=(\rpt{n-1}{t-1},1)$ or  
  $\underline{\lambda}=(3,\rpt{n-3}{t-1})$.
\end{Proposition}

\medskip
\begin{proof}
It is clear that 
\begin{itemize}
\item if the last part of $\rpt{n-1}{t-1}$ is not equal to $1$, then $\DD{\rpt{n-1}{t-1},1}=(n-1,1)$ and 
\item if the first part of $\rpt{n-3}{t-1}$ is at most $2$, then $\DD{3,\rpt{n-3}{t-1}}=(n-1,1)$.
\end{itemize}
 Suppose now, $\underline{\lambda}=(\lambda_1,\lambda_2,\ldots,\lambda_t)$
 and $\DD{\underline{\lambda}}=(n-1,1)$. By Lemma \ref{thm:lemma1}, we have that 
 $\underline{\lambda}=(\lambda_1,\lambda_2,\ldots,\lambda_s,\lambda_{s+1},\ldots,\lambda_t)$,
where 
$\lambda_1+\lambda_2+\ldots+\lambda_s=n-1$ or $2s+\lambda_{s+1}+\ldots+\lambda_t=n-1$.
In the first case, clearly, $\underline{\lambda}=(\rpt{n-1}{s},1)$. In the second case we have that
$\sum_{i=1}^{s}(\lambda_i-2)=1$ and thus, $s=1$ and $\lambda_1=3$.
\end{proof}

\bigskip

\begin{Proposition}\label{thm:-3}
If $2\leq r\leq 5$   and $\mu-r \geq 1$, then
 $$\left|\DDi{\mu,\mu-r}\right|=(r-1)(\mu-r).$$

Moreover, 
$$\DDi{\mu,\mu-2}=\left\{ (\mu,\lambda_2,\lambda_3,\ldots,\lambda_t); \; 
  (\lambda_2,\lambda_3,\ldots,\lambda_t) \, \text{ is almost rectangular}, \, \mu-\lambda_t \geq 2\right\},$$

 $$\DDi{\mu,\mu-3}=\left\{ (\mu-\varepsilon,\lambda_2,\lambda_3,\ldots,\lambda_t); \; 
 \lambda_2-\lambda_t\leq 1, \, \mu-\varepsilon-\lambda_t \geq 2, \, \varepsilon \in \{0,1\} \right\}$$
and 
  \begin{eqnarray*}
 \DDi{\mu,\mu-4}&=&\left\{ (\mu_1,\mu_2,\lambda_3,\lambda_4,\ldots,\lambda_t); \; 
    \mu_1-\mu_2 \leq 1, \, \lambda_3-\lambda_t\leq 1, \, \mu_1-\lambda_t \geq 2 \right\} \cup \\
 && \cup \; \left\{ (\mu-\varepsilon,\lambda_2,\lambda_3,\ldots,\lambda_t); \; \lambda_2-\lambda_t\leq 1, \, 
    \mu-\varepsilon-\lambda_t \geq 2, \, \varepsilon \in \{0,2\} \right\}.
 \end{eqnarray*}
 \end{Proposition}
 
\medskip
\begin{proof}
As a corollary of Theorem \ref{thm:dif2} we obtain the set 
$\DDi{\mu,\mu-2}$  and $\left|\DDi{\mu,\mu-2}\right|=\mu-2$.

 If $\DD{\underline{\lambda}}=(\mu,\mu-3)$,  then $\underline{\lambda}$ is of the form of Lemma \ref{thm:lemma1}, where $s=1$. 
 By \eqref{eq:12}, we get
 $\lambda_1 = 2\mu-3-(\lambda_2+\lambda_3+\ldots+\lambda_t) \geq 2\mu-3-(\mu-2)=\mu-1$ and from \eqref{eq:11} it follows that 
 $\lambda_1 \leq \mu$. Thus, we have that either $\lambda_1=\mu$ or $\lambda_1=\mu-1$. 
 It can be easily verified that all such partitions $\underline{\lambda}$ satisfy the condition 
 $\DD{\underline{\lambda}}=(\mu,\mu-3)$.
 
 If $\lambda_1=\mu$, then $(\lambda_2,\lambda_3,\ldots,\lambda_t)$ is an arbitrary almost rectangular partition of
 $\mu-3$. If $\lambda_1=\mu-1$, then $(\lambda_2,\lambda_3,\ldots,\lambda_t)$ is an almost rectangular partition of
 $\mu-2$, where $t \geq 3$. (Note, that otherwise $\lambda_1=\lambda_t+1$.) Therefore, $\left|\DDi{\mu,\mu-3}\right|=2(\mu-3)$.

\smallskip

 If $\DD{\underline{\lambda}}=(\mu,\mu-4)$, then by Lemma \ref{thm:lemma1},  it follows that $s \leq 2$. 
 
 If $s=1$, then by \eqref{eq:12} we have that
 $\lambda_1 = 2\mu-4-(\lambda_2+\lambda_3+\ldots+\lambda_t) \geq 2\mu-4-(\mu-2)=\mu-2$.
 Since by \eqref{eq:11}, $\lambda_1 \leq \mu$ we consider 3 cases. 
 
 If $\lambda_1=\mu$, then
 $(\lambda_2,\lambda_3,\ldots,\lambda_t)$ is an almost rectangular partition of $\mu-4$ and by \cite[Thm. 7]{KoOb09}, 
 it follows that $\DD{\mu,\lambda_2,\lambda_3,\ldots,\lambda_t}=(\mu,\mu-4)$. 
 If $\lambda_1=\mu-1$ and 
 $(\lambda_2,\lambda_3,\ldots,\lambda_t)$ is an almost rectangular partition of $\mu-3$,
 then   $\DD{\mu-1,\lambda_2,\lambda_3,\ldots,\lambda_t}=(\mu-1,\mu-3)$ and thus
 no such partition is in $\DDi{\mu,\mu-4}$. 
 If $\lambda_1=\mu-2$, then
 $(\lambda_2,\lambda_3,\ldots,\lambda_t)$ is an almost rectangular partition of $\mu-2$. If, in addition, $t \geq 3$, then
 $\DD{(\mu-2,\lambda_2,\lambda_3,\ldots,\lambda_t)}=(\mu,\mu-4)$.
 
  If $s=2$, then, by \eqref{eq:12}, we have that
  $\lambda_1+\lambda_2 = 2\mu-4-(\lambda_3+\lambda_4+\ldots+\lambda_t) \geq 2\mu-4-(\mu-4)=\mu$.
  Thus, $(\lambda_1,\lambda_2)$ is an almost rectangular partition of $\mu$ and $(\lambda_3,\lambda_4,\ldots,\lambda_t)$
  is an almost rectangular partition of $\mu-4$, such that $\mu_1-\lambda_t \geq 2$. This is true if and only if $t \geq 4$.
  
  Now, it is easy to compute that $\left|\DDi{\mu,\mu-4}\right|=\mu-4+\mu-3+\mu-5=3(\mu-4)$.
\end{proof}

\bigskip

\begin{Question}
Is it true that
  $$\left|\DDi{\mu,\mu-r}\right|=(r-1)(\mu-r)$$
for all $r \geq 5$?
\end{Question}

\bigskip

One can also ask a question, what are maximal and minimal partitions in $\DDi{\underline{\mu}}$? 
Clearly, the maximal partition in $\DDi{\underline{\mu}}$ is $\underline{\mu}$. However, there is not a unique minimal partition
in $\DDi{\underline{\mu}}$, as shown in the next example.

\bigskip

\begin{Example}
One can easily check that 
 $$\DDi{6,2}=\{(6,2), \, (6,1^2), \, (4,2^2), \, (4,2,1^2), \, (4,1^4),\, (3^2,1^2)\}$$
and that there are 2 minimal partitions $(3^2,1^2)$ and $(4,1^4)$ in $\DDi{6,2}$.
\end{Example}

\bigskip

Recall that the \emph{rank} of partition $(\lambda_1,\lambda_2,\ldots,\lambda_t) \in \partition{n}$
is defined as the number $n-t$.
So, the partition with the  minimal rank is the partition with the most parts. Now, we can prove the following:

\bigskip

\begin{Proposition}
 For every $r \geq 2$, the partition $(\mu+2,1^{\mu+r-2})$ is in $\DDi{\mu+r,\mu}$ and this is the unique partition 
 with the  minimal rank in $\DDi{\mu+r,\mu}$.
\end{Proposition}

\smallskip

\begin{proof}
 Since $r \geq 2$, we have by \cite[Thm. 7]{KoOb09}, that $\DD{\mu+2,1^{\mu+r-2}}=(\mu+r,\mu)$.
 
 Suppose that $\DD{\underline{\lambda}}=(\mu+r,\mu)$ and that $\underline{\lambda}$ has a rank at most 
 $\mu+1$, which is the rank of $(\mu+2,1^{\mu+r-2})$. Then, $\underline{\lambda}$ is of the form
 $(\lambda_1,\lambda_2,\ldots\lambda_s,\lambda_{s+1},\ldots,\lambda_t)$,
 where $t \geq \mu+r-1$, $\lambda_1-\lambda_s\leq 1$, $\lambda_{s+1}-\lambda_t\leq 1$ and $\lambda_1- \lambda_t \geq 2$.
 Now, we have
 \begin{eqnarray*}
  \lambda_1+\lambda_2+\ldots+\lambda_s &\leq& \mu+r \\
  2s+\lambda_{s+1}+\lambda_{s+2}+\ldots+\lambda_t &\leq& \mu+r\, .
 \end{eqnarray*} 
 If $t \geq \mu+r$, it follows that $\mu+r-s \leq t-s \leq \lambda_{s+1}+\lambda_{s+2}+\ldots+\lambda_t \leq \mu+r-2s$,
 which is a contradiction.
 Otherwise, if $t=\mu+r-1$, then $\mu+r-1-s = t-s \leq \lambda_{s+1}+\lambda_{s+2}+\ldots+\lambda_t \leq \mu+r-2s$ and
 thus $s=1$. Now, we conclude that $\lambda_2=\lambda_3=\dots=\lambda_{\mu+r}=1$ and 
 $\underline{\lambda}=(\mu+2,1^{\mu+r-2})$.
 Thus, $(\mu+2,1^{\mu+r-2})$ is the unique partition with rank equal to $\mu+r-1$ and no partition in $\DDi{\mu+r,\mu}$ has
 greater rank.
\end{proof}

\bigskip

\begin{Question}
 Let $(\mu_1,\mu_2,\ldots,\mu_s)$ be a stable partition. Is it true that the partition with the minimal rank in
 $\DDi{\mu_1,\mu_2,\ldots,\mu_s}$ is equal to $(\mu_2+2,\mu_3+2,\ldots,\mu_s+2,1^{\mu_1-2(s-1)})$?
\end{Question}

\bigskip
\bigskip


\end{document}